# On Double-Entry Bookkeeping: The Mathematical Treatment

David Ellerman

University of California/Riverside

Home page: www.ellerman.org

Email: david@ellerman.org

**Abstract**

Double-entry bookkeeping (DEB) implicitly uses a specific mathematical construction, the group of differences using pairs of unsigned numbers ("T-accounts"). That construction was only formulated abstractly in mathematics in the 19th century—even though DEB had been used in the business world for over five centuries. Yet the connection between DEB and the group of differences (here called the "Pacioli group") is still largely unknown both in mathematics and accounting. The precise mathematical treatment of DEB allows clarity on certain conceptual questions and it immediately yields the generalization of the double-entry method to multi-dimensional vectors typically representing the different types of property involved in an enterprise or household.

**Introduction**

Double-entry bookkeeping (DEB) was developed during the fifteenth century and was published in 1494 as a system by the Italian mathematician Luca Pacioli (Pacioli 1914).[1] Double-entry bookkeeping has been used for over five centuries in commercial accounting systems. If the mathematical formulation of any field should be well understood, one would think it might be accounting. Remarkably, however, the mathematical formulation of double entry accounting—algebraic operations on ordered pairs of numbers—is largely unknown in mathematics as well as in accounting.

DEB implicitly uses a very specific and precise mathematical construction (that is now part of undergraduate abstract algebra). That construction was developed in the nineteenth century by William Rowan Hamilton as a result of using ordered pairs to deal with the complex numbers (Hamilton 1837). The multiplicative version of this construction is the "group of fractions" which uses ordered pairs of whole numbers (written vertically) to enlarge the system of positive whole numbers to the system of positive fractions containing multiplicative inverses (just reverse the entries in a fraction to get its inverse). The ordered pairs construction that is relevant to DEB is the additive case usually called the "group of differences."[2] It is used to construct a number system with additive inverses by using operations on ordered pairs of non-negative or, to be more precise, unsigned numbers. The key to grasp the connection with DEB is to make the identification:

$$\text{ordered pairs of numbers in group of differences construction}$$
$$= \text{two-sided T-accounts of DEB.}$$

In view of this identification, the group of differences (or fractions in the multiplicative case) might be called the *Pacioli group* (Ellerman 1982), or if we "violate" Stigler's Law of eponymy,[3] the *Pacioli-Cotrugli group*.

In spite of some attention to DEB by mathematicians [e.g., (DeMorgan 1869), (Cayley 1894), and (Kemeny, Schleifer, Snell, *et al.* 1962)], this connection has *not* been noted in mathematics (not to mention in accounting) with one *perhaps solitary* exception. In a semi-popular book, D. E. Littlewood noted the connection:

> The bank associates two totals with each customer's account, the total of moneys credited and the total of moneys withdrawn. The net balance is then regarded as the same if, for example, the credit amounts of £102 and the debit £100, as if the credit were £52 and the debit £50. If the debit exceeds the credit the balance is negative.
>
> This model is adopted in the definition of signed integers. Consider pairs of cardinal numbers (a, b) in which the first number corresponds to the debit, and the second to the credit. (Littlewood 1960 p. 18)

---

[1] It should be noted that the double-entry system seems to be first described, but not first published, by the Dubrovnik (then Ragusa) merchant-economist Benedetto Cotrugli (or Benedikt Kotruljevic in the Croatian version) in 1458 (Yamey 1994). Pacioli was familiar with and credited Cotrugli 's manuscript. "With a fine sense of ethics, Paciolo credits Cotrugli (1458) with having originated the method of double entry bookkeeping." (Langer 1958 p. 483)

[2] See (Jacobson 1951 p. 10); (Littlewood 1960 p. 18); (Dubisch 1965 p. 17); (Mac Lane and Birkhoff 1967 p. 45); and so forth.

[3] "No scientific discovery is named after its original discoverer." (Stigler 1999 p. 277)

Arthur Cayley (1821-1895) came close. In the year before his death, he published a small pamphlet *The Principles of Book-keeping by Double Entry* in which he wrote:

> The Principles of Book-keeping by Double Entry constitute a theory which is mathematically by no means uninteresting: it is in fact like Euclid's theory of ratios an absolutely perfect one, and it is only its extreme simplicity which prevents it from being as interesting as it would otherwise be. (Cayley 1894)

In the pamphlet, Cayley only described double-entry bookkeeping in the practical informal terms familiar from his fourteen years of work as a lawyer. In his Presidential Address to the British Association for Advancement of Science, he hinted that the "notion of a negative magnitude" is "used in a very refined manner in bookkeeping by double entry" (Cayley 1896 p. 434). In neither place did Cayley present the ordered pairs treatment of T-accounts. However the reference to the theory of ratios in Euclid points to exactly that connection, i.e., to the group of multiplicative ordered pairs called "fractions" or "ratios."

> If four numbers be proportional, the number produced from the first and fourth will be equal to the number produced from the second and third; [and conversely]. (Euclid 1926)

The translator and commentator, Sir Thomas Heath, points out that, in modern notation, this is "If a:b = c:d, then ad = bc; and conversely" [Commentary in (Euclid 1926 p. 319)]. This is precisely the equality condition (equality of cross-multiples) in the multiplicative group of fractions (see Table 1).

Except for Littlewood's remark (which was not further developed to describe the double-entry method) and Cayley's pointer, the author has not been able to find a single mathematics (not to mention, accounting) book or paper, elementary or advanced, popular or esoteric, which [independently of (Ellerman 1982), (Ellerman 1985), (Ellerman 1986)] shows how the ordered pairs of the group of differences construction are the T-accounts used in the business world for about five centuries.

**The Pacioli Group**

Double-entry accounting is based on the group of differences or Pacioli group construction using unsigned numbers.[4] It generalizes naturally to multi-dimensional double-entry accounting where the unsigned or non-negative numbers are replaced by vectors of non-negative numbers (i.e., ordered lists of non-negative numbers). The objects in the group of differences are ordered pairs of unsigned numbers which are identified with the *T-accounts* of DEB. The left-hand side (LHS) number d is the debit entry and the right-hand side (RHS) number c is the credit entry.

*T-account: [ d // c ] = [ debit number // credit number ].*[5]

---

[4] The Pacioli group and its use to carry out the operations of DEB is a rational reconstruction, not an effort at historical exegesis of Pacioli or Cotrugli (where the author has no expertise). For instance, it may well be that Pacioli implicitly assumed a complete set of accounts (so the double-entry principle would always hold) without seeing those accounts as coming from an equation such as today's balance sheet equation (Sangster 2010).

[5] The double-slash notation was suggested by Pacioli. "At the beginning of each entry, we always provide 'per', because, first, the debtor must be given, and immediately after the creditor, the one separated from the other by two little slanting parallels (virgolette), thus, //,… ." (Pacioli 1914 p. 43)

The algebraic operations on T-accounts are much like the operations on the ordered pairs of numbers called "fractions" except that addition is substituted for multiplication. In order to illustrate the additive-multiplicative analogy between T-accounts and fractions, the basic definitions will be developed in parallel columns. For a fraction or multiplicative T-account,[6] we may take the numerator as the debit entry and the denominator as the credit entry.

| **Table 1** | **Additive Case** | **Multiplicative Case** |
|---|---|---|
| Operation on T-accounts | T-accounts add together by adding debits to debits and credits to credits [ x // y ] + [ w // z ] = [ x+w // y+z ]. | Fractions multiply together by multiplying numerator times numerator and denominator times denominator (x/y)(w/z) = (xw/yz). |
| Identity element | The identity element for addition is the zero T-account [ 0 // 0 ]. | The identity element for multiplication is the unit fraction (1/1). |
| Equality between two T-accounts. | Given two T-accounts [ x // y ] and [ w // z ], the *cross-sums* are the two numbers obtained by adding the credit entry in one T-account to the debit entry in the other T-account. The equivalence relation between T-accounts is defined by setting two T-accounts *equal* if their cross-sums are equal: [ x // y ] = [ w // z ]  if x+z = y+w. | Given two fractions (x/y) and (w/z), the *cross-multiples* are the two numbers obtained by multiplying the numerator of one with the denominator of the other. The equivalence relation between fractions is defined by setting two fractions *equal* if their cross-multiples are equal: (x/y) = (w/z) if xz = yw. |
| Inverses | The negative or additive inverse of a T-account is obtained by reversing the debit and credit entries: – [ x // y ] = [ y // x ]. | The multiplicative inverse of a fraction is obtained by reversing the numerator and denominator: $(x/y)^{-1}$ = (y/x). |
| "Disjointness" of two T-accounts | Given two non-negative numbers x and y, let min(x,y) be the minimum of the two numbers. Two non-negative numbers x and y are said to be *disjoint* if min(x,y) = 0. | Given two integers w and x, let gcd(w,x) be the greatest common divisor of w and x (the largest integer dividing both). Two integers w and x are said to be *relatively prime* if gcd(w,x) = 1. |
| "Reduced form" | A T-account [x // y] is in *reduced form* if x and y are disjoint. | A fraction (x/y) is in *lowest terms* if x and y are relatively prime. |
| Unique reduced form representation | Every T-account [x // y] has a unique reduced form representation [x–min(x,y) // y–min(x,y)]. | Each fraction (x/y) has a unique representation in lowest terms (x/gcd(x,y))/(y/gcd(x,y)). |
| Example | Consider the T-account [12 // 5]. The minimum of the debit and credit numbers is 5 so the reduced form representation is [7 // 0]. | Consider the fraction (28/35). The greatest common divisor of the numerator and denominator is 7 so the fraction in lowest terms is (4/5). |

The *group of differences* or *Pacioli group* **P** = **P**(ℕ) of the natural numbers ℕ = {0,1,2,…} consists of the ordered pairs [x // y], called *T-accounts* or *T-terms*, of unsigned whole numbers from ℕ with the above definition of addition and equality.[7] The Pacioli group **P** = **P**(ℕ) is isomorphic to (i.e., is in one-to-one correspondence in an addition-preserving way with) the integers $\mathbb{Z} = \{\ldots,-2, -1, 0, 1, 2, \ldots\}$, $\mathbf{P}(\mathbb{N}) \cong \mathbb{Z}$. But it is important to realize that it is isomorphic to the integers in *two* different ways. In the integers $\mathbb{Z}$, we can distinguish the positive and negative whole numbers but it is only multiplication that distinguishes the two. In the Pacioli group $\mathbf{P}$, there is only addition. There are no negative numbers, only unsigned numbers in T-terms.[8] However there are additive inverses, for instance [x // y] and [y // x] are additive inverses since they add to the:

$$[x // y] + [y // x] = [0 // 0].$$
*Zero T-account*

---

[6] There is in fact a whole system of multiplicative DEB (Ellerman 1982 pp. 58–66) to update multiplicative equations of natural numbers like 216×140 = 12×63×40.

[7] We could work with real or rational numbers but for simplicity will stick to the whole numbers.

However, this does *not* mean that either [x // y] or [y // x] is a "negative number." For instance, neither [7 // 5] nor [5 // 7] should be identified as the negatively signed number –2, and, for that matter, neither should be identified with the positively signed number +2. It depends in the context of the double-entry method on whether the T-account is debit-balance or credit-balance. In mathematical terms, there are *two* isomorphisms with associate T-accounts [x // y] with signed numbers. There is the *debit isomorphism*, which associates [x // y] with x–y, and the *credit isomorphism*, which associates [x // y] with y–x.[9] It is often said that DEB avoids negative numbers, but, to be precise, it also avoids positive numbers; it uses unsigned numbers.

In order to translate from T-accounts [x // y] back and forth to the signed integers in $\mathbb{Z}$, one needs to specify whether to use the debit or credit isomorphism. This will be done in DEB by labeling each balance sheet (and income statement) T-account as *debit-balance* or *credit-balance*. Thus if a T-account [x // y] is debit-balance, the corresponding number or "balance" in the account is x–y, and if it is credit-balance, then the corresponding balance is the number y–x. All these calculations with numbers or "scalars"[10] will generalize below to vectors or ordered lists of unsigned numbers to give the multi-dimensional vector version of DEB.

**The Double-Entry Method: Scalar Case**

**What is the "double" in double-entry bookkeeping?**

Since one of the clarifications resulting from the mathematical formulation of DEB is understanding just what is "double" that is unique to DEB. Hence we first visit the standard treatment of that in the accounting texts. Given an equation x + ... + y = w + ... + z, it is not possible to change just one term in the equation and have it still hold. Two or more terms must be changed. The fact that two or more terms (or "accounts") must be changed is actually *not* (NB: *not*) unique to the double-entry method even though that is a commonplace in the accounting literature.

A business transaction does not affect just one item alone. There are at least two items to be considered in each transaction....The dual aspect of each transaction forms the basis underlying

---

[8] There is a whole literature about DEB and negative numbers. While Pacioli as a mathematician knew about negative numbers (Scorgie 1989), it is uncontroversial to say that the use of DEB by businesses in that time was facilitated by not requiring familiarity with the relative new notion of a negative number (Heeffer 2011).

[9] In going from x and y in $\mathbb{N}$ to x–y or y–x in $\mathbb{Z}$, we are, as usual, identifying the unsigned numbers in $\mathbb{N}$ with the non-negative numbers in $\mathbb{Z}$.

[10] In the context of vectors, single numbers are called "scalars."

what is called *double-entry accounting*. (Moore and Jaedicke, 1967, pp. 648–649)

> Every event that is recorded in the accounts affects at least two items; there is no conceivable way of making only a single change in the accounts. Accounting is therefore properly called a "double-entry" system. (Anthony 1970 p. 32)

That mathematical fact (of affecting two or more terms) is a characteristic of the transaction itself (the changes in the equation), not a characteristic of the method of recording the transaction. The double-entry method is a method of encoding an equation using ordered pairs of unsigned numbers, i.e., T-accounts, to record transactions and make changes in the equation. The characteristic doubleness of the double-entry method is the *double-sidedness* of the T-accounts and the mathematical properties that follow (e.g., equal debits and credits in a transaction, and equal debits and credits in the trial balance of the whole set of accounts or Ledger).

To show that the doubleness of entries is not unique to DEB, we need to precisely give the alternative *system* of recording transactions. The alternative to the double-entry method is to record a transaction by making an entry of adding a *signed* (positive or negative) number to each affected (single-sided) account. The exact same two or more accounts in the equation would still always be affected by this alternative method of recording a transaction (since that is a property of the transaction itself, not a property of the recording method). Such a system using a complete balance sheet of accounts is a *complete* accounting system to record and update the balance sheet equation. Since this complete single-sided account method of recording transactions has:

- no two-sided T-accounts (only single-sided accounts),
- no debits or credits of unsigned numbers (using signed numbers instead),
- no double entry principle in the form of equal debits and credits in a transaction (although the same mathematics appears in a different form), and
- no trial balance by adding debits and credits,

it is *not* the double-entry system (if language is to be used with any precision at all).

This alternative system using single-sided accounts and "double-sided" or signed numbers might be called "single-entry accounting" in analogy with the double-sided accounts of double-entry accounting. Unfortunately, the phrase "single-entry accounting" is usually used to denote simply an *incomplete* accounting "system," e.g., small businesses that don't keep track of any equity account so some transactions will only affect a single account (since the equity account is ignored). But without an equation, that is not an alternative "system" at all. The usual accounting-text juxtaposition of a complete double-entry system with an incomplete "single-entry system" is a false comparison. The real comparison is between:

- using unsigned ("single-sided") numbers in double-sided accounts (i.e., DEB) or
- using signed ("double-sided") numbers in "single-sided" accounts (i.e., the complete system without debits and credits or the other characteristic features of DEB).

For brevity, we might refer to that alternative system as the *SSS system* (or *SSS method*) for Single-Sided accounts with Signed numbers (from $\mathbb{Z}$) in contrast to the *DSU system* (double-entry method) using Double-Sided accounts (in the Pacioli group **P**) with Unsigned numbers (from $\mathbb{N}$).

The point is that the "double" of the double-entry system is typically identified with the fact that two or more accounts are affected when that *same* characteristic of the transaction (as opposed to a characteristic of the recording method) *also* holds for the SSS system of single-sided accounts and signed numbers.

To make matters perhaps even more confusing to the layperson, the two systems are equivalent in the sense that given the same beginning equation and the same transactions, both systems will give the same ending equation. Otherwise one system would just be wrong instead of being an alternative recording and updating method. The virtues of DEB (as opposed to the SSS system) lie not in its mathematical properties but in its human properties: dealing with unsigned numbers (e.g., the unintuitiveness of taking either asset or equity accounts as negative—even after negative numbers were well-known), double-entry rule of equal debits and credits for transactions, and trial balance of equal debits and credits to check equations. Computers operate with signed numbers (since they have to multiply numbers too) so an accounting program might well have a double-entry human interface using debits and credits—which are then internally mapped (using either the debit or the credit isomorphism) to signed numbers and single-sided accounts (i.e., single memory locations).

We now turn to the mathematical structure of DEB (using scalars). The same transactions will be illustrate later using the SSS system, and then both the DEB (or DSU) and SSS systems will be generalized to multi-dimensional vectors.

**The mathematics of double-entry bookkeeping**

The Pacioli group is only the "algebra" where double-entry method does its calculations; it is not the double-entry method itself. The double-entry method is a procedure (using the Pacioli group explicitly or implicitly) to start with an equation (additive in the case of DEB), record changes in the terms of the equation, and correctly obtain the new ending equation. In the case of DEB, the equation is the balance sheet equation and the changes in the equation result from transactions affecting the balance sheet accounts (and income statement accounts if those are used).

Consider an example of a company with the simplified initial balance sheet equation:

$$\text{Assets} = \text{Liabilities} + \text{Equity}$$
$$15000 = 10000 + 5000.$$

*Beginning Scalar Balance Sheet*

It is customary in accounting (although not mathematically necessary) to move each term or "account" to the side of the equation so that it is preceded by a plus sign. A T-account equal to the zero T-account $[0 \mathbin{/\mkern-5mu/} 0]$ is called a *zero-account* or *zero-term*. Equations, like the equation of unsigned numbers $x = y + z$, encode as zero-accounts. Each left-hand side (LHS) term x is encoded as a debit-balance T-account $[x \mathbin{/\mkern-5mu/} 0]$ and each right-hand side (RHS) term y is encoded as a credit-balance T-account $[0 \mathbin{/\mkern-5mu/} y]$. These T-accounts then would add up to the zero-account $[0 \mathbin{/\mkern-5mu/} 0]$:

$$[x \mathbin{/\mkern-5mu/} 0] + [0 \mathbin{/\mkern-5mu/} y] + [0 \mathbin{/\mkern-5mu/} z] = [x \mathbin{/\mkern-5mu/} y+z] = [0 \mathbin{/\mkern-5mu/} 0]$$

*An equational zero-account*

(where the last equation between T-terms holds because of the equality of cross-sums). The balance sheet equation thus encodes as an equation zero-account. The *Ledger* in DEB is just the listing of the T-accounts of the balance sheet zero-account leaving out the plus signs.[11]

|  | Assets | Liabilities | Equity |
|---|---|---|---|
|  | [15000 // 0] | [0 // 10000] | [0 // 5000] |

*Beginning Ledger of T-Accounts*

Consider three transactions in a productive firm.

1. $1200 of input inventories are used up and charged directly to equity.
2. $1500 of product is produced, sold, and added directly to equity.
3. $800 principal payment is made on a loan.

Each transaction is then encoded as a *transactional zero-term* and added to the appropriate terms of the equational zero-account. For instance, the first transaction subtracts 1200 from Assets and subtracts 1200 from Equity. The Assets account is encoded as a LHS or debit-balance account so the subtracting of a number from it would be encoded as adding the T-account [0 // 1200] to it. Equity is encoded as a RHS or credit-balance term so subtracting 1200 from it would be encoded as adding [1200 // 0] to it. The other transactions are encoded in a similar manner. The list of the transaction zero-terms is just the *Journal* and the encoding of the transactions as zero-terms is the double-entry principle (i.e., each transaction must be recorded with equal debits and credits).

|  | Assets | Liabilities | Equity |
|---|---|---|---|
| Original equation zero-account: | [15000 // 0] | [0 // 10000] | [0 // 5000] |
| +Transaction 1 zero-term: | [0 // 1200] |  | [1200 // 0] |
| +Transaction 2 zero-term: | [1500 // 0] |  | [0 // 1500] |
| +Transaction 3 zero-term: | [0 // 800] | [800 // 0] |  |
| = Ending equation zero-account: | [16500 // 2000] | [800 // 10000] | [1200 // 6500] |
| = (in reduced form) | [14500 // 0] | [0 // 9200] | [0 // 5300]. |

*Initial Ledger + Journal = Ending Ledger*

The initial T-accounts in the Ledger add up to the zero-account (initial trial balance). Each transaction is encoded as two or more T-accounts that add to the zero-account (double entry principle). Zero added to zero equals zero and, since zero is the additive identity, *only* zero added to zero gives zero.[12] Thus adding the transaction zero-terms to the initial equation zero-account (posting Journal to the Ledger) will yield another equation zero-account (which can be checked by taking another trial balance) representing the ending balance sheet.[13] Each T-account is then decoded according to whether it was encoded as debit-balance or credit-balance to obtain the ending balance sheet equation in its usual form.

---

[11] For simplicity and brevity, we are leaving out the temporary or flow accounts of the income statement and dealing directly with Equity—rather than accumulating the changes to Equity in the income statement and then closing it into Equity at the end of the accounting period. To use income statement accounts, one would simply add a credit-balance Revenue account and a debit-balance Expenses account to the Ledger (with zero initial balances). The two transactions 1 and 2 affecting Equity would then be posted in the usual way to the Expense and Revenue accounts. Finally another transaction or two at the end of the period would "close the income statement into the balance sheet" by making the changes to Equity and resetting the zero balance in the Revenue and Expense accounts.

[12] Zero being the additive identity means for any x, 0 + x = x + 0 = x. Hence if 0 + x = 0 then x = 0.

[13] We have left out the income statement accounts for simplicity but if they are used, then those temporary accounts are closed into the balance sheet by the usual closing transactions to get the ending balance sheet.

$$\begin{array}{ccccc} \text{Assets} & = & \text{Liabilities} & + & \text{Equity} \\ 14500 & = & 9200 & + & 5300. \end{array}$$

*Ending Balance Sheet Equation*

**Transaction "Matrices"**

The mathematical structure of the Pacioli group of T-accounts (using scalars) can be presented using the visual aid of an M×M square table (where there are M different accounts) often called a *transactions matrix*.[14] The transactions table was first used by Augustus DeMorgan who called it a "table of *double-entry*" (DeMorgan 1869 p. 83) but was popularized a century later by the mathematicians John Kemeny and colleagues (Kemeny, Schleifer, Snell, *et al.* 1962) (who were unaware of the Pacioli group treatment of DEB). The two sides of transaction T-terms for, say, Assets are then visually presented with the credit side in the *column* labeled "Assets" and the debit side in the *row* labeled "Assets."

In our simplified example, the transactional zero-term:

$$[0 // 1200] + [1200 // 0]$$
$$\quad \textit{Assets} \qquad \textit{Equity}$$

would be represented in a transaction table by rotating the Assets T-term with the credit entry of 1200 by 90 degrees so that it overlaps with the equal debit entry of 1200 in the Equity T-term.

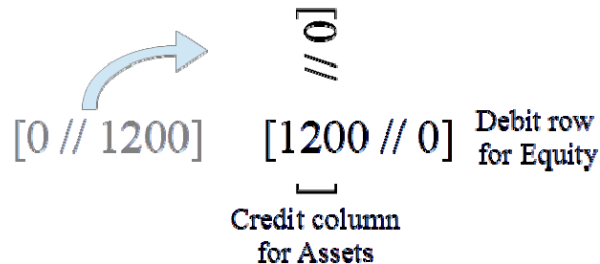

Figure 1: *Rotating the credits to represent a transaction zero-term as a table entry*

Thus the doubly-classified unsigned number 1200 is just a visually clever way to represent that transaction zero-term [0 // 1200] as the change in Assets and [1200 // 0] as the change in Equity.

Here are the three transactions from the previous example.

| Table 2 | Assets (Cr.) | Liabilities (Cr.) | Equity (Cr.) |
|---|---|---|---|
| Assets (Dr.) | | | 1500 |
| Liabilities (Dr.) | 800 | | |
| Equity (Dr.) | 1200 | | |

*Transactions table*

---

[14] Technically a "matrix" is a way of presenting a linear transformation from one space of vectors to another so the key operation is the *multiplication* of a matrix times a vector to yield another vector. The so-called "transaction matrices" used in the table-presentation of double-entry bookkeeping *do not multiply* so it may be impressive but ultimately misleading to call the tables "matrices." There are genuine uses of matrices in business mathematics; see (Kemeny, Schleifer, Snell, *et al.* 1962) or (Shank, 1972) for applications to proportional overhead allocation and Markov chains.

The transactions table does not avoid using the T-terms but only represents them in a doubly-classified way. When the T-terms are used explicitly, one finds out how much to change the balance in the accounts by adding up the Journal zero-terms:

|  | Assets | Liabilities | Equity |
|---|---|---|---|
| Transaction 1 zero-term: | [0 // 1200] |  | [1200 // 0] |
| +Transaction 2 zero-term: | [1500 // 0] |  | [0 // 1500] |
| +Transaction 3 zero-term: | [0 // 800] | [800 // 0] |  |
| = Total transaction zero-term | [1500 // 2000] | [800 // 0] | [1200 // 1500] |

*Sum of Journal transactions*

In the transactions table approach, one implicitly does exactly the same T-term additions but with a different visual presentation. In the above sum of Journal transactions, there are three sums of the right-hand side credit entries and those are precisely the three sums of the credit columns in the transaction table. And the three sums of the left-hand side debit entries are precisely the three sums of the debit rows in the transaction Table 2.

But what do we do with, say, the two numbers 1500 and 2000 obtained by summing the column and row for Assets?

In DEB using explicit T-accounts, the transactions sum is added to the beginning Assets account [15000 //0]:

$$\begin{array}{c} \text{Assets} \\ [15000 // 0] \\ + [1500 // 2000] \\ \hline = [16500 // 2000] \\ = [14500 // 0] \end{array}$$

*Updating Assets account*

(where the last equation is by the equality of cross-sums) which gives the ending Asset account as a T-account. Since Assets is a debit-balance account, this is decoded as an ending balance of +14500 in the Assets account.

In the transactions table presentation, we have for, say, the Assets account, the change in debits of 1500 and the change in credits of 2000. This is just an implicit T-term [1500 // 2000] for the changes in the account and we have to know the debit- or credit-balance nature of the account to know which way to net the two numbers. For instance, if the $i^{th}$ account is debit-balance then since:

> the entries in the $i^{th}$ row of the ledger matrix represent all the debits to account number i and entries in the $i^{th}$ column represent all the credits to this account, the net impact of transactions during a period on account number i can be computed by subtracting the $i^{th}$ column sum from the $i^{th}$ row sum. (Shank 1972 p. 22)

That is, for a debit-balance T-account, we apply the debit isomorphism where [x // y] gives x−y. If the $i^{th}$ account was a credit-balance account, then the two unsigned numbers would be subtracted in the opposite way. That is, for a credit-balance T-account, we apply the credit isomorphism where [x // y] gives y−x. Since the Assets account is debit-balance, the column-sum credit of 2000 is subtracted from the row-sum debit of 1500 to get the same change in the Assets account balance: 15000 + (1500 − 2000) = 14500. Similarly for the other accounts.

The point is that the transaction table approach is making the exact same calculations using the exact same information about debits and credits as the explicit DEB approach, but it is all presented in a different visual manner and with the background information about accounts being debit- or credit-balance being used at the crucial moment to decide what to do with the unsigned column-sum and row-sum for an account.

This implicit use of the T-account algebra is not a different mathematical approach from the double-entry method using explicit T-accounts; it's just implicit. Furthermore, the transactions table presentation obscures the underlying action of the Pacioli group by keeping its use implicit, so the treatments in the literature [e.g., (Kemeny, Schleifer, Snell, *et al.* 1962) or (Shank 1972)] do not connect it explicitly to the group of differences construction (Pacioli group) or give the then-obvious generalization to multi-dimensional accounting.[15]

**The mathematics of the SSS system of bookkeeping**

In view of the isomorphism(s) between the T-accounts of unsigned numbers in the Pacioli group $\mathbf{P}(\mathbb{N})$ and the signed numbers in the integers $\mathbb{Z} = \{\ldots,-2,-1,0,1,2,\ldots\}$, there are two ways to construct a corresponding SSS system by applying either the credit or debit isomorphisms. If we indulge the layperson's understanding that "to debit" means "to subtract" then we would choose the credit isomorphism so that [x // y] maps to y–x. Or we can indulge the idea that all asset accounts should be seen as positive and apply the debit isomorphism so that [x // y] maps to x–y.

I will use the debit isomorphism since that is, in effect, used in the best mathematical treatment of accounting in the recent literature by Robert Nehmer and colleagues (Cruz Rambaud, Perez, Nehmer, *et al.* 2010; Nehmer and Robinson 1997). The previous example of the mathematical treatment of DEB is then transformed into the SSS system by applying the debit isomorphism to each equation and transaction zero-term.

For instance, the beginning balance sheet zero-account: [15000 // 0] [0 // 10000] [0 // 5000], transforms by the debit isomorphism to:

| Assets | Liabilities | Equity |
|--------|-------------|--------|
| 15000  | –10000      | –5000  |

*Beginning Ledger of Single-Sided accounts and Signed numbers*

(where it would be excessively pedantic to write a positive signed number like 1500 as +1500). The fact that the balance sheet was represented by a zero T-account now appears as the fact that the row (or column in Nehmer's presentation) sums to zero. Call such a row a *zero-row*. In the same manner, the transactional zero-terms are transformed into rows of signed numbers that sum to zero (which is the form of the double-entry principle of recording a transaction with equal debits and credits on the other side of the isomorphism) and are added to the beginning balance sheet accounts to get the ending accounts.

---

[15] The table or "matrix" presentation of the scalar case of DEB does not generalize well to the multi-dimensional vector case as is explained in (Ellerman 1982, pp. 132–134).

|                              | Assets | Liabilities | Equity |
|------------------------------|--------|-------------|--------|
| Original equation zero-row:  | 15000  | –10000      | –5000  |
| +Transaction 1 zero-row:     | –1200  |             | +1200  |
| +Transaction 2 zero-row:     | +1500  |             | –1500  |
| +Transaction 3 zero-row:     | –800   | +800        |        |
| = Ending equation zero-row:  | 14500  | –9200       | –5300  |

*Initial Ledger + Journal = Ending Ledger*
*using Single-Sided accounts and Signed numbers*

The result is the ending balance sheet in the form: Assets–Liabilities–Equity = 0 which can easily be put in the usual form: A = L + E.

For the purposes of the "doubleness" question, this accounting system using single-sided accounts (and "double-sided" or signed numbers) has exactly the same characteristic of each transaction affecting two or more accounts, so that cannot be the *unique* characteristic of the DEB (or DSU) system using double-sided accounts (and unsigned numbers).

Since the two systems are equivalent as systems for recording transactions to update equations, it is hard to resist swinging back and forth across the (debit) isomorphism to "enrich" (or muddle) the SSS system with "debit" and "credit" terminology—as Nehmer and company unfortunately do.[16] Thus I have gone to some pains to be conceptually clear about the two distinct systems while being equally clear about the isomorphism(s) between them.[17]

**The Double-Entry Method: the General Multi-dimensional Case**

**Previous failed attempts at multi-dimensional DEB**

One acid test of a mathematical formulation of a theory is whether or not it facilitates the generalization of the theory. Normal bookkeeping does not deal with incommensurate physical quantities; everything is expressed in the common units of money (which is often taken as a necessary condition for the double-entry method). Is there a generalization of DEB to deal with multi-dimensional physical incommensurates such as the different types of goods and services without collapsing the "apples and oranges" to one quantity of value?[18] If successful, this multi-dimensional DEB would show, contrary to many texts, that the common monetary unit is not a prerequisite for DEB.

---

[16] In spite of being the best treatment of mathematical accounting in the recent literature, there are some infelicities in the treatment of these questions by Nehmer and colleagues (Cruz Rambaud, Perez, Nehmer, *et al.* 2010) in addition to casually going back and forth across the debit isomorphism (which muddles the SSS and DSU systems). For instance, they use single-sided accounts but then define "T-diagrams" (2010 p. 53)—which are then used only to represent the transactions and not the accounts themselves. They do consider (2010 pp. 54–55) the Pacioli group from the author's earlier work (1986, 1985, 1982). But even that definition is somewhat infelicitous since they allow all the negative numbers of the integers $\mathbb{Z}$ to appear in the T-accounts of the Pacioli group which negates the whole point of the group of differences construction. That is, instead of constructing $\mathbf{P}(\mathbb{N}) \cong \mathbb{Z}$ from $\mathbb{N}$, they construct $\mathbf{P}(\mathbb{Z}) \cong \mathbb{Z}$ from $\mathbb{Z}$. More importantly, they consider only the Pacioli group and ignore the whole double-entry method of T-accounts with the balance sheet, income statement, Ledger, Journal, posting the Journal to the Ledger, and so forth—not to mention the multi-dimensional generalization. Hence it seems the treatment of DEB using the group of differences construction is *still* essentially unknown in the accounting literature.

[17] The "affects-two-accounts" mischaracterization of double-entry bookkeeping is far too embedded in the DNA of the accounting profession to realistically expect it to be changed now. At best, on occasions where there is clarity in distinguishing the DEB (or DSU) system from the SSS system, then retronymic phrases will probably be used so that the DEB system will be called the "double-entry system with T-accounts" while the SSS system will be called the "double-entry system with single-sided accounts" or the "double-entry system without T-accounts."

In the literature on the mathematics of accounting, there was a proposed "solution" to this question, a system of multi-dimensional physical accounting (Ijiri 1967, 1966, 1965). In this system, most of the normal structure of DEB was lost:

- there was no balance sheet equation,
- there were no equity or proprietorship accounts,
- the temporary or nominal accounts could not be closed, and
- the "trial balance" did not balance.

It is common for certain aspects of a theory to be lost in a generalization of the theory. The accounting community had apparently accepted the failure of all these features of DEB as the necessary price to be paid to generalize DEB to incommensurate physical quantities. For example, the systems of "Double-entry multidimensional accounting" previously published in the accounting literature [see also (Charnes, Colantoni and Cooper 1976), or (Haseman and Whinston 1976)] had acquiesced in the absence of the balance-sheet equation.

> For instance, the convenient idea of an accounting identity is lost since the dimensional and metric comparability it assumes is no longer present except under special circumstances. (Charnes, Colantoni and Cooper 1976 p. 333)

Yet when DEB is mathematically formulated using the group of differences, then the generalization to vectors of incommensurate physical quantities is immediate and trivial. *All* of those normal features of DEB—such as the balance-sheet equation, the equity account, the temporary accounts, and the trial balance—are preserved in the generalization as indicated in the multi-dimensional example given here (Ellerman 1986, 1985, 1982). Thus the "accepted" generalized model of DEB was simply a failed attempt at generalization which had been "received" as a successful generalization that unfortunately had to "sacrifice" certain features of DEB.

**The Pacioli group using vectors**

The general double-entry method starts with an equation between sums of n-dimensional vectors. Vector equations such as:

$$(6, -3, 10) + (-2, 5, -2) = (4, 2, 8)$$

hold if the sums of the corresponding components of the vectors give scalar equations, i.e.,

$$6-2 = 4, -3+5 = 2, \text{ and } 10-2 = 8.[19]$$

As in the scalar case, double-entry vector accounting is a way of recording changes in equations due to "transactions" changing the various physical quantities. There is, as before, a way to do this using single-sided accounts and vectors with signed components. But our goal is to first generalize the double-entry method where there are double-sided T-accounts using vectors with unsigned components.[20]

---

[18] By "multi-dimensional DEB" I do not mean using scalar DEB coupled with "multi-dimensional" database records; instead I mean a multi-dimensional version of the full DEB systems with multi-dimensional (vector) equations and with transactions to update the equations. Hence all the usual machinery of balance sheets, T-accounts, Ledgers, Journals, flow-statements, and so forth will generalize. Moreover, the property vector accounting system, in a certain sense, underlies the value scalar accounting system.

[19] The numbers will later be interpreted as types of property (goods and services) but for now the focus is just on the mathematics.

Vector equations are first *encoded* in the Pacioli group constructed from ordered pairs of n-dimensional vectors with unsigned components (instead of the unsigned scalars used in scalar DEB). If we stick to whole numbers as before, then the n-dimensional vectors would be from the n-fold product: $\mathbb{N}^n = \mathbb{N} \times \ldots \times \mathbb{N}$ of all ordered n-tuples $x = (x_1, \ldots, x_n)$ for unsigned whole numbers $x_i$ from $\mathbb{N}$. The Pacioli group $\mathbf{P} = \mathbf{P}(\mathbb{N}^n)$ of $\mathbb{N}^n$ is constructed exactly as before with the definitions in Table 1 but the x, y, w, and z are now vectors in $\mathbb{N}^n$. Then instead of the two isomorphisms between $P(\mathbb{N})$ and $\mathbb{Z}$, we have the two (debit and credit) isomorphisms between $\mathbf{P}(\mathbb{N}^n)$ and $\mathbb{Z}^n = \mathbb{Z} \times \ldots \times \mathbb{Z}$.

A vector is *non-negative* if all its components are non-negative (the unsigned numbers in $\mathbb{N}$ are traditionally identified with the non-negative numbers in $\mathbb{Z}$). Since the entries in the vectors in a T-account must be unsigned (or non-negative), we must first develop a way to separate out the positive and negative components of a vector from $\mathbb{Z}^n$.

Given two vectors $x = (x_1, \ldots, x_n)$ and $y = (y_1, \ldots, y_n)$, let max(x,y) be the vector with the maximum of $x_i$ and $y_i$ as its $i^{th}$ component, and let min(x,y) be the vector with the minimum of $x_i$ and $y_i$ as its $i^{th}$ component. Two non-negative vectors x and y are said to be *disjoint* if min(x,y) = 0. The *positive part* of a vector x is non-negative vector $x^+ = \max(x,0)$, the maximum of x and the zero vector [note that "0" is used, depending on the context, to refer to the zero scalar or the zero vector]. For $x = (6, -3, 10)$, $x^+ = (6,0,10)$. The *negative part* of x is the non-negative vector $x^- = -\min(x,0)$, the negative of the minimum of x and the zero vector so that $x^- = (0,3,0)$. The positive and negative parts of a vector x are disjoint non-negative vectors.[21] Every vector x has a *Jordan decomposition* $x = x^+ - x^-$:

$$x = (6, -3, 10) = (6, 0, 10) - (0, 3, 0) = x^+ - x^-.$$
*Example of Jordan decomposition of a vector*

The debit and credit isomorphisms were previously presented as going from the Pacioli group $P(\mathbb{N})$ to $\mathbb{Z}$ using the mappings of [x // y] to x−y or to y−x respectively. The vector versions use the exact same formulas to go from the Pacioli group $\mathbf{P}(\mathbb{N}^n)$ to $\mathbb{Z}^n$. But being isomorphisms, there is also the two isomorphisms going in the opposite direction from $\mathbb{Z}^n$ to $\mathbf{P}(\mathbb{N}^n)$. The two isomorphisms that map $\mathbb{Z}^n$ vectors to T-accounts of $\mathbb{N}^n$ vectors are the:

*debit isomorphism*: x in $\mathbb{Z}^n$ maps to $[x^+ // x^-]$ in $\mathbf{P}(\mathbb{N}^n)$, and

*credit isomorphism*: x in $\mathbb{Z}^n$ maps to $[x^- // x^+]$ in $\mathbf{P}(\mathbb{N}^n)$.

---

[20] Expanding the transaction amount in scalar accounting in a link to database information such as the account affected, the date, or the counterparty, and then calling the result a "vector" does not create "vector accounting" (Mattessich 1964 p. 94) in the sense used here. Vectors, unlike database records, can be meaningfully added together.

[21] An "abuse of language" is involved in calling the non-negative vector $x^-$ the "negative part" of x.

Two vector T-accounts [x // y] and [w // z] are *equal* if (as with scalars) their cross-sums are equal:

$$[x // y] = [w // z] \text{ if and only if } x+z = y+w.$$
*Definition of equality of vector T-accounts*

A vector T-account [x // y] is said to be in *reduced form* if x and y are disjoint (i.e., min(x,y) = 0). Every vector T-account [x // y] is equal to a vector T-account in reduced form, namely

$$[x // y] = [x–\min(x,y) // y–\min(x,y)].$$
*Putting a vector T-account into reduced form*

(as one can see by computing the cross-sums). Given a vector equation, $x + ... + y = w + ... + z$, each left-hand side (LHS) vector x is encoded via the debit isomorphism as a debit-balance T-account $[x^+ // x^-]$ and each right-hand side (RHS) vector w is encoded via the credit isomorphism as a credit-balance T-account $[w^- // w^+]$. Then the original equation holds if and only if the sum of the encoded T-accounts is the vector zero-account [0 // 0] so it is called an *equation zero-account*:

$$x + ... + y = w + ... + z$$
if and only if
$$[x^+ // x^-] + ... + [y^+ // y^-] + [w^- // w^+] + ... + [z^- // z^+] = [0 // 0].$$
*Encoding an Equation as an Equation Zero-Account*

This could be illustrated by applying it to the vector equation:

$$(6, –3, 10) + (–2, 5, –2) = (4, 2, 8).$$

The Jordan decompositions are;

    $x = (6, –3, 10)$ so $x^+ = (6,0,10)$ and $x^- = (0,3,0)$;

    $y = (–2, 5, –2)$ so $y^+ = (0,5,0)$ and $y^- = (2,0,2)$; and

    $w = (4,2,8)$ so $w^+ = (4,2,8)$ and $w^- = (0,0,0)$.

Hence the equation zero-account is:

    $[(6,0,10) // (0,3,0)] + [(0,5,0) // (2,0,2)] + [(0,0,0) // (4,2,8)] = [(0,0,0) // (0,0,0)].$

To check the equation, one adds up the debit vectors on the LHS to get;

$$(6,0,10) + (0,5,0) + (0,0,0) = (6,5,10),$$

and adds up the credit vectors on the RHS to get:

$$(0,3,0) + (2,0,2) + (4,2,8) = (6,5,10)$$

so the cross-sums are both $(6,5,10) + (0,0,0) = (6,5,10)$ and thus are equal.

Since only plus signs can appear between the T-accounts in an equational zero-account, the plus signs can be left implicit. The listing of the T-accounts in an equational zero-account (without the plus signs) is the *Ledger*.

Changes in the various accounts in the beginning equation are recorded as *transactions*. Transactions must be recorded as valid algebraic operations which transform equations into equations. Since equations encode as zero-accounts, a valid algebraic operation would transform zero-accounts into zero-accounts. There is only one such operation in the Pacioli group: add on a zero-account or zero-term. Zero plus zero equals zero. The zero-accounts representing transactions are called *transaction zero-terms*. The listing of the transactional zero-terms is the *Journal*.

A series of valid additive operations on a vector equation can then be presented in the following standard scheme:

    Beginning Equation Zero-Account
+   <u>Transaction Zero-Terms</u>
=   Ending Equation Zero-Account

or, in more conventional terminology,

   Beginning Ledger
+ <u>Journal</u>
= Ending Ledger.

The process of adding the transaction zero-terms to the initial Ledger to obtain the Ledger at the end of the accounting period is called *posting the Journal to the Ledger*. The fact that a transaction zero-term is equal to [0 // 0] is the *double-entry principle* that transactions are recorded with equal debits and credits. The summing of the debit and credit sides of what should be an equation zero-account to check that it is indeed a zero-account is the *trial balance*. All those features from scalar case of DEB carry over to the general vector case.

At the end of the cycle, the ending equational zero-account is put into reduced form and decoded to obtain the equation that results from the algebraic operations represented in the transactions. T-accounts [x // y] that are debit-balance are decoded as x–y on the left side of the equation, and T-accounts [x // y] that are credit-balance are decoded as y–x on the right side of the equation. In an accounting application, the T-accounts in the final equation zero-account would be partitioned as debit-balance and credit-balance according to the side of the initial equation from which they were encoded.

The mathematics is general but in the accounting context, the n-dimensional vectors would typically represent the different types of property rights. Indeed, the Pacioli group formulation of DEB was originally developed to show that the stocks and flows of different types of property in a firm could be described using the double-entry method (Ellerman 1982). We will consider a simple model where there are only three types of property: cash, outputs (widgets), and inputs (half-widgets). These goods will be listed in that order in each three-dimensional vector.

Let the initial asset vector be (9000, 40, 50) so the firm has $9000 cash, 40 units of widgets in the output inventory, and 50 units of half-widgets in the input inventory. The firm also has a $10000 liability represented by the vector (10000, 0, 0) so the equity vector (Assets – Liabilities) is given by the (net) property vector (–1000, 40, 50). Thus the initial balance sheet (vector) equation is:

        Assets          =  Liabilities       +  Equity
       (9000, 40, 50)    =  (10000, 0, 0)    +  (–1000, 40, 50).

*Initial Vector Balance-Sheet Equation*

This is encoded as the following equation zero-account or vector Ledger:

    Assets                    Liabilities              Equity
  [(9000, 40, 50)//(0, 0, 0)]    [(0, 0, 0)//(10000, 0, 0)]    [(1000, 0, 0)//(0, 40, 50)].

*Initial Vector T-Accounts in Ledger*

The underlying production process is very simple. Two units of the half-widget inputs are combined to make one widget. Hence the following physical transactions underlie the value transactions recorded in scalar accounting (where we split the production and sale of the outputs as the transactions 2a and 2b).

1. 30 units of the half-widgets inputs are used up in production.
2a. 15 units of the widgets output are produced.
2b. 15 units of the widgets output are sold for $100 each.
3. $800 principal payment is made on a loan.

These transactions are then encoded as transaction zero-terms and added to the Ledger T-accounts. For instance, the using-up of 30 units of half-widgets input is recorded as crediting 30 input units to Assets and debiting the 30 units to Equity (we again deal, for simplicity, directly with Equity rather than using "income statement" or property-flow accounts) to give the transaction zero-term:

$$[(0, 0, 0)//(0, 0, 30)] + [(0, 0, 30)//(0, 0, 0)] = [(0,0,0)//(0,0,0)].$$
*Assets*           *Equity*

The other transactions are encoded similarly to form the vector-accounting version of the Journal (listing of transaction zero-terms). Then the Journal is added ("posted") as usual to the Ledger representing the beginning equation zero-account to yield the ending equation zero-account.

| | Assets | Liabilities | Equity |
|---|---|---|---|
| | [(9000, 40, 50)//(0, 0, 0)] | [(0, 0, 0)//(10000, 0, 0)] | [(1000, 0, 0)//(0, 40, 50)] |
| 1. | [(0, 0, 0)//(0, 0, 30)] | | [(0, 0, 30)//(0, 0, 0)] |
| 2a. | [(0, 15, 0)//(0, 0, 0)] | | [(0, 0, 0)//(0, 15, 0)] |
| 2b. | [(1500, 0, 0)//(0, 15, 0)] | | [(0, 15, 0)//(1500, 0, 0)] |
| 3. | [(0, 0, 0)//(800, 0, 0)] | [(800, 0, 0)//(0, 0, 0)] | |
| | ================ | ================ | ================ |
| | [(10500, 55, 50)//(800, 15, 30)] | [(800, 0, 0)//(10000, 0, 0)] | [(1000, 15, 30)//(1500, 55, 50)] |
| = | [(9700, 40, 20)//(0, 0, 0)] | [(0, 0, 0)//(9200, 0, 0)] | [(0, 0, 0)//(500, 40, 20)] |

*Initial Vector Ledger + Vector Journal = Ending Vector Ledger*

where the last line of Ledger accounts is in reduced form. The reduced accounts are then decoded to obtain the ending balance-sheet equation:

| Assets | | Liabilities | | Equity |
|---|---|---|---|---|
| (9700, 40, 20) | = | (9200, 0, 0) | + | (500, 40, 20). |

*Ending Vector Balance-Sheet Equation*

Given a set of prices (or valuation coefficients), the vectors can be evaluated so that the vector accounts of property accounting collapse to the scalar accounts of value accounting. Thus the property vector accounting may be said to *underlie* the usual value accounting. For instance, suppose that the prices per unit (cash, output, input) are (1, 100, 40) respectively. Multiplying the physical quantities times their price and adding up (called the "scalar product" or "dot product" in mathematics) yields the value of the property vectors. For instance, the ending Assets property vector (9700, 40, 20) has the value:

$$(1, 100, 40) \bullet (9700, 40, 20) = 1 \times 9700 + 100 \times 40 + 40 \times 20 = 9700 + 4000 + 800 = 14500.$$

*Scalar product of a price vector times a property vector*

In this manner, we obtain the balance-sheet equation in the previous example of value or scalar accounting.

$$\begin{array}{ccccc} \text{Assets} & = & \text{Liabilities} & + & \text{Equity} \\ 14500 & = & 9200 & + & 5300. \end{array}$$

*Scalar Equation = Price Vector times Vector Equation*

Thus we see how property accounting can use double-entry accounting with vectors to trace out the property transactions that underlie the value transactions recorded in conventional accounting.

**Vector accounting with Single-Sided accounts and vectors of Signed numbers**

As in the scalar case, the mathematical treatment of DEB is easily converted to the corresponding SSS system by applying one of the isomorphisms. Applying the debit isomorphism, which takes [x // y] to x–y (as before), gives the following:

| Assets | Liabilities | Equity |
|---|---|---|
| (9000, 40, 50) | (–10000, 0, 0) | (1000, –40, –50) |

*Initial Vector Single-Sided Accounts Ledger*

which is a zero-row of vectors (i.e., the vectors add to zero) so it is a vector equation in the form: A–L–E = 0. Applying the debit isomorphism to the Journal (transaction rows also adding to zero) and posting it to the Ledger gives the vector version of those operations using the SSS system.

| | Assets | Liabilities | Equity |
|---|---|---|---|
| | (9000, 40, 50) | (–10000, 0, 0) | (1000, –40, –50) |
| 1. | (0, 0, –30) | | (0, 0, 30) |
| 2a. | (0, 15, 0) | | (0, –15, 0) |
| 2b. | (1500, –15, 0) | | (–1500, 15, 0) |
| 3. | (–800, 0, 0) | (800, 0, 0) | |
| | ============ | ============ | ============ |
| | (9700, 40, 20) | (–9200, 0, 0) | (–500, –40, –20). |

*Initial Vector Ledger + Vector Journal = Ending Vector Ledger*
*using Single-Sided accounts of vectors with Signed numbers*

That is the ending balance sheet vector equation in the A–L–E = 0 form. The use of vectors (of signed numbers) to represent the stocks and flows of property rights (goods and services) with evaluations by scalar products with price vectors is completely standard in mathematical economics [see any text such as (Varian 1984)] which includes but is not restricted to input-output theory.

**Concluding Remarks**

There is a precise mathematical system underlying the system of double-entry bookkeeping used in businesses for over five centuries. The underlying mathematical construction (the group of differences) was only explicitly developed in mathematics itself in the $19^{th}$ century. Today the mathematical treatment of DEB [in spite of being published a third of a century ago (Ellerman 1982)], is little known in both the fields of mathematics and accounting—which shows a remarkable insulation between the fields.

The precise mathematical treatment of DEB helps to clear up some old canards in the accounting literature by showing that what is specifically "double" in DEB is the two-sidedness of the T-accounts (which gives the system its debits and credits), not the fact that every transaction affects two or more accounts. It also resolves the question of whether or not a common unit of account (e.g., money) is necessary to use the double-entry method since the mathematics extends easily to a vector version of DEB (e.g., representing the different types of goods and services involved in business transactions) where no common measure of value is assumed.